 \documentclass[11pt]{article}
 \usepackage[top=1in,bottom=1in,left=1.5in,right=1.5in]{geometry}
  \usepackage{amsmath,amssymb}
  \usepackage{latexsym}
 \usepackage[dvips]{pstricks} 
  \usepackage{pst-node}
  \usepackage[dvips]{graphicx}

  \newcommand{\C}{\mathbb{C}}
  
  \newcommand{\N}{\mathbb{N}}
  
  \newcommand{\R}{\mathbb{R}}
  \newcommand{\Z}{\mathbf{Z}}

  \newcommand{\hs}{\hspace*{\parindent}}

  \newcommand{\tr}{\mathop{\mathrm{tr}}\nolimits}

  \newcommand{\qed}{\hspace*{\fill} $\Box$\\}

  \newtheorem{theo}{\bfseries \hs Theorem}[section]
  \newtheorem{defn}[theo]{\bfseries \hs Definition}

  \numberwithin{equation}{section} 

\begin{document}

 \title{A note on the nonzero spectrum of irreducible matrices}
 \author{
 Shmuel Friedland\\
 Department of Mathematics, Statistics and Computer Science\\
 University of Illinois at Chicago\\ Chicago, Illinois 60607-7045,
 USA\\ \texttt{E-mail: friedlan@uic.edu}
 }

 \renewcommand{\thefootnote}{\arabic{footnote}}
 \date{March 19, 2011 }
 \maketitle
 \begin{abstract}
 In this note we extend the necessary and sufficient conditions of Boyle-Handleman  1991 and Kim-Ormes-Roush 2000 for a nonzero eigenvalue
 multiset of primitive matrices over $\R_+$ and $\Z_+$, respectively, to irreducible matrices.
 \end{abstract}



 \section{Introduction}\label{intro}
 Denote by $\R^{n\times n}\supset \R_+^{n\times n}$ the algebra of real valued $n\times n$ matrices and the cone of $n\times n$ nonnegative matrices,
 respectively.  For $A\in\R^{n\times n}$ denote by $\Lambda(A)=\{\lambda_1(A),\ldots,\lambda_n(A)\}$ the eigenvalue multiset of $A$, i.e.
 $\det (zI-A)=\prod_{i=1}^n (z-\lambda_i(A))$.  An outstanding problem in matrix theory, called NIEP, is to characterize a
 multiset $\Lambda=\{\lambda_1,
 \ldots,\lambda_n\}$ which is an eigenvalue multiset of some $A\in \R^{n\times n}$.  Denote
 by $\rho(\Lambda):=\max\{|\lambda|, \lambda\in\Lambda\}$, and by $\Lambda(r)$ all elements in $\Lambda$ satisfying $|\lambda|=r\ge 0$.
 For $\lambda\in\Lambda$ denote by $m(\lambda)\in\N$ the multiplicity of $\lambda$ in $\Lambda$.
 The obvious necessary conditions for $\Lambda=\Lambda(A)$ for some $A\in\R^{n\times n}$ are the trace conditions:
 \begin{equation}\label{tracecon}
 s_k(\Lambda):=\sum_{i=1}^n \lambda_i^k \ge 0 \textrm{ for } k=1,\ldots ,
 \end{equation}
 since $s_k(\Lambda(A))=\tr A^k$.
 The following theorem is deduced straightforward from \cite[Thm 2]{Fri78}.  (See \S2.)
 \begin{theo}\label{fritrcon}.  Let $\Lambda=\{\lambda_1,\ldots,\lambda_n\}$ be a multiset of complex numbers.
 Assume that the inequalities in (\ref{tracecon}) hold except for a finite number values of $k$.
 Then
 \begin{enumerate}
 \item\label{fritrcon1} $\bar\Lambda=\Lambda$.
 \item\label{fritrcon2}
 $\rho(\Lambda)\in \Lambda$.
 \item\label{fritrcon3}
 $m(\rho(\Lambda))\ge m(\lambda)$
 for all $\lambda\in \Lambda(\rho(\Lambda))$.
 \item\label{fritrcon4}
 Assume that $\rho(\Lambda)>0$ and let $\{\lambda_1,\ldots,\lambda_p\}$ be all
 distinct elements of $\Lambda$ such that $|\lambda_i|=\rho(\Lambda)$ and $m(\lambda_i)=m(\rho(\Lambda))$
 for $i=1,\ldots,p$.  Then $\zeta\Lambda(\rho(\Lambda))=\Lambda(\rho(\Lambda))$ for $\zeta=e^{\frac{2\pi\sqrt{-1}}{p}}$.
 \end{enumerate}
 \end{theo}

 By considering the diagonal elements of $A^k$ and comparing them with the diagonal elements of $A^{km}$
 Loewy and London added additional the following necessary conditions \cite{LoL78}.
 \begin{theo}\label{loelcond}  Let $\Lambda=\{\lambda_1,\ldots,\lambda_n\}$ be an eigenvalue multiset of some
 $A\in\R_+^{n\times n}$.  Then in addition to the inequalities (\ref{tracecon}) the following inequalities hold.
 \begin{equation}\label{loelcond1}
 s_{km}(\Lambda)\ge \frac{1}{n^{k-1}}(s_{m}(\Lambda))^k \textrm{ for } m,k-1=1,\dots.
 \end{equation}
 In particular,
 \begin{equation}\label{loelcond2}
 \textrm{if } s_{m}(\Lambda) >0 \textrm{ then } s_{km}(\Lambda) >0 \textrm{ for } k=2,\ldots .
 \end{equation}
 \end{theo}

 The inequalities (\ref{tracecon}) and (\ref{loelcond1}) imply that $\Lambda$ is an eigenvalue multiset of some
 $A\in\R_+^{n\times n}$ in the following cases: $n=3$; $n=4$ and $\Lambda$ is a multiset of real numbers.
 For $n=4$ and nonreal $\Lambda=\{\lambda_1,\lambda_2,\lambda_3,\lambda_4\}$ the conditions (\ref{tracecon}) and (\ref{loelcond1})
 are not sufficient \cite{LoL78}.  The necessary and sufficient conditions are given in \cite{TM07}.
 The inequality $ns_4(\Lambda)\ge (s_2(\Lambda))^2$ in (\ref{loelcond1})
 can be improved to $(n-1)s_4(\Lambda)\ge (s_2(\Lambda))^2$ if $s_1(\Lambda)=0$ and $n$ odd \cite{LaM98}.
 \begin{defn}\label{defnfrobset}  A multiset $\Lambda=\{\lambda_1,\ldots,\lambda_n\}\subset \C$
 is called a Frobenius multiset
 if the following conditions hold.
 \begin{enumerate}
 \item\label{defnfrobset2}  $\bar \Lambda=\Lambda$.
 \item\label{defnfrobset3}  $\rho(\Lambda)\in \Lambda$.
 \item\label{defnfrobset4}  $m(\lambda)=1$ for each $\lambda\in \Lambda(\rho(\Lambda))$.
 \item\label{defnfrobset5}  Assume that $\#\Lambda(\rho(\Lambda))=p$.
 Then $\zeta \Lambda=\Lambda$ for $\zeta=e^{\frac{2\pi\sqrt{-1}}{p}}$.
 \end{enumerate}
 \end{defn}
 The Frobenius theorem for irreducible $A\in\R_+^{n\times n}$, i.e. $(I+A)^{n-1}$ is a positive matrix,
 claims that $\rho(\Lambda(A))>0$ and $\Lambda(A)$ is a Frobenius set.  In particular, an irreducible $A\in\R_+^{n\times n}$
 is primitive, i.e. $A^{(n-1)^2+1}$ is a positive matrix, if and only if $\Lambda(\rho(\Lambda))=\{\rho(\Lambda)\}$,
 see \cite[XII.\S5]{Gan59} and \cite[\S8.5.9]{HoJ88}.

 We say that a multiset $\Lambda=\{\lambda_1,\ldots,\lambda_n\}$, where $\lambda_i\ne 0$ for $i=1,\ldots,n$,
 is a nonzero eigenvalue multiset of a nonnegative matrix if there exist an integer $N\ge n$ and $A\in\R_+^{N\times N}$,
 such that $\Lambda$ is obtained from $\Lambda(A)$ be removing all zero eigenvalues.  The following remarkable
 theorem was proved by Boyle and Handelman \cite{BoH91}.  Namely, a multiset $\Lambda\subset\C\backslash\{0\}$ is a nonzero spectrum of a nonnegative
 primitive matrix if and only if $\Lambda(\rho(\Lambda))=\{\rho(\Lambda)\}$, and the inequalities (\ref{tracecon}) and
 (\ref{loelcond2}) hold.  See the recent proof of Thomas Laffey \cite{Laf10} of a simplified version of this result.
 The aim of this note to extend the theorem of Boyle-Handelman to a nonzero eigenvalue multiset of
 nonnegative irreducible matrices.
 \begin{theo}\label{nscnzspirmat}  Let $\Lambda$ be a multiset of nonzero complex numbers.  Then $\Lambda$ is a nonzero eigenvalue
 multiset of a nonnegative irreducible matrix if and only if $\Lambda$ is a Frobenius set, and
 (\ref{tracecon}) and (\ref{loelcond2}) hold.
 \end{theo}
 Similarly, we extend the results of Kim, Ormes and Roush \cite{KOR00} to a nonzero eigenvalue multiset of nonnegative irreducible matrices
 with integer entries.
 \section{Proofs of Theorems \ref{fritrcon} and \ref{nscnzspirmat}}
 \textbf{Proof of Theorem \ref{fritrcon}.}  For $\Lambda=\{0,\ldots,0\}$ the theorem is trivial.
 Assume that $\rho(\Lambda)>0$.  Consider the function
 $$f_{\Lambda}(z)=\sum_{i=1}^n \frac{1}{1-\lambda_i z}=\sum_{k=0}^{\infty} s_k(\Lambda)z^k.$$
 Assume that $s_k(\Lambda)\ge 0$ for $k> N$.  Then by subtracting a polynomial $P(z)$ of degree $N$ at most,
 we deduce that $f_0(z):=f(z)-P(z)$ has real nonnegative MacLaurin coefficients.  So $\overline{f_0(\bar z)}=f(z)$.
 Hence $\bar\Lambda=\Lambda$.
 The radius of convergence of this series is $R(f_{\Lambda})=\frac{1}{\rho(\Lambda)}$.  The principal part of $f$ is
 $f_1:=\sum_{i, |\lambda_i|=\rho(\Lambda)} \frac{1}{1-\lambda_i z}$.  So $\pi(f_{\Lambda}(z))=\{(\lambda_1,m(\lambda_1),1),\ldots,
 (\lambda_q,m(\lambda_q),1)\}$, where $\lambda_1,\ldots,\lambda_q$ are all pairwise distinct elements of $\Lambda(\rho(\Lambda))$, see
 \cite[Dfn. 1]{Fri78}.  Then parts \ref{fritrcon2}--\ref{fritrcon4} follow from \cite[Thm. 2]{Fri78}.\\

 \noindent
 \textbf{Proof of Theorem \ref{nscnzspirmat}.}  Assume first that $\Lambda$ is a nonzero eigenvalue multiset of a nonnegative irreducible
 matrix.  The Frobenius theorem yields that $\Lambda$ has to be a Frobenius set, and (\ref{tracecon}) and (\ref{loelcond2}) hold.
 Assume now that $\Lambda$ has is a Frobenius set, and (\ref{tracecon}) and (\ref{loelcond2}) hold.
 In view of the Boyle-Handelman theorem it is enough to consider the case
 \begin{equation}\label{nontrcase}
 \Lambda(\rho(\Lambda))=\{\rho(\Lambda),\zeta \rho(\Lambda),\ldots,\zeta^{p-1} \rho(\Lambda)\}, \textrm{ for }
 \zeta=e^{\frac{2\pi\sqrt{-1}}{p}} \textrm{ and } 1<p\in\N.
 \end{equation}

 Observe first that $s_k(\Lambda)=0$ if $p\not | k$.  Let $\phi:\C\to\C$ be the map $z\mapsto z^p$.  Since $\zeta \Lambda=\Lambda$
 it follows that for $z\in\Lambda$ with multiplicity $m(z)$ we obtain the multiplicity $z^p$ in $\phi(\Lambda)$ is $p m(z)$.
 Hence $\phi(\Lambda)$ is a union of $p$ copies of a Frobenius set $\Lambda_1$, where $\rho(\Lambda_1)=\rho(\Lambda)^p$ and
 $\Lambda_1(\rho(\Lambda_1))=\{\rho(\Lambda_1)\}$.
 Moreover $s_{kp}(\Lambda)=ps_k(\Lambda_1)$.   Hence $\Lambda_1$ satisfies the assumptions of the Boyle-Handelman theorem.
 Thus there exists a primitive matrix $B\in\R_+^{M\times M}$ whose nonzero eigenvalue multiset is $\Lambda_1$.    Let $A=[A_{ij}]_{i=j=1}^p$
 be the following nonnegative matrix of order $pM$.
 \begin{equation}\label{defirrA}
 A=\left[\begin{array}{cccccc} 0_{n\times n}& I_n&0_{n\times n}&0_{n\times n}&\ldots&0_{n\times n}\\
 0_{n\times n}&0_{n\times n}& I_n&0_{n\times n}&\ldots&0_{n\times n}\\
 \vdots&\vdots&\vdots&\vdots&\vdots&\vdots\\
 0_{n\times n}&0_{n\times n}&0_{n\times n}&0_{n\times n}&\ldots&I_n\\
 B&0_{n\times n}&0_{n\times n}&0_{n\times n}&\ldots&0_{n\times n}
 \end{array}\right].
 \end{equation}
 Then $A$ is irreducible and the nonzero part of eigenvalue multiset $\Lambda(A)$ is $\Lambda$. \qed
 \section{An extension of Kim-Ormes-Roush theorem}
 In this section we give necessary and sufficient conditions on a multiset $\Lambda$ of nonzero complex number to be a
 nonzero eigenvalue multiset of a nonnegative irreducible matrix with integer entries.
 Recall the M\"obius function $\mu:\N\to\{-1,0,1\}$.  First $\mu(1)=1$.  Assume that $n>1$. If $n$ is not square free,
 i.e. $n$ is divisible by $l^2 $ for some positive integer $l>1$, then $\mu(n)=0$.  If $n>1$ is square free, let $\omega(n)$
 be the number of distinct primes that divide $n$.  Then $\mu(n)=(-1)^{\omega(n)}$.  The following theorem is a generalization
 of the Kim-Ormes-Roush theorem \cite{KOR00}.
 \begin{theo}\label{KORgen}  Let $\Lambda$ be a multiset of nonzero complex numbers.  Then $\Lambda$ is a nonzero eigenvalue
 multiset of a nonnegative irreducible matrix with integer entries if and only if the following conditions hold.
 \begin{enumerate}
 \item\label{KORgen1}   $\Lambda$ is a Frobenius set.
 \item\label{KORgen2}  The coefficients of the polynomial $\prod_{\lambda\in\Lambda}(z-\lambda)$ are integers.
 \item\label{KORgen3} $t_k(\Lambda):=\sum_{d|k}\mu(\frac{k}{d}) s_d(\Lambda)\ge 0$ for all $k\in\N$.

 \end{enumerate}

 \end{theo}
 The case $\Lambda(\rho(\Lambda))=\{\rho(\Lambda)\}$ is the Kim-Ormes-Roush theorem.\\

 \noindent
 \textbf{Proof.}  Assume that $\Lambda$ is a nonzero spectrum of a nonnegative irreducible matrix with integer entries,
 i.e. $A\in\Z_+^{N\times N}$.
 Then part \ref{KORgen1} follows from the Frobenius theorem.  Since $\det(zI-A)$ has integer coefficients we deduce we deduce
 part \ref{KORgen2}.  It is known that $t_k(\Lambda)=t_k(\Lambda(A))$ is the number of minimal loops of length $k$ in the
 directed multigraph induced by $A$, see \cite{BoH91}.  Hence part \ref{KORgen3} holds.

 Suppose that $\Lambda$ satisfies \ref{KORgen1}--\ref{KORgen3}.  In view of the Kim-Ormes-Roush theorem it is enough to assume the case
 (\ref{nontrcase}).  We now use the notations and the arguments of the proof of Theorem \ref{nscnzspirmat}.   First $s_k(\Lambda)=0$ if $p\not|k$.
 Second $\prod_{\lambda\in\Lambda}(z-\lambda)=\prod_{\kappa\in\Lambda_1}(z^p-\kappa)$.  Hence $\prod_{\kappa\in\Lambda_1}(z-\kappa)$
 has integer coefficients.
 A straightforward calculation shows that $t_{pk}(\Lambda)=pt_{k}(\Lambda_1)$.  Hence $t_{k}(\Lambda_1)\ge 0$.
 Kim-Ormes-Roush theorem yields the existence of $B\in\Z_+^{M\times M}$ such that $\Lambda_1$ is the nonzero eigenvalue
 multiset of $B$.  Hence $\Lambda$ is the nonzero eigenvalue set of $A\in \Z_+^{pM\times pM}$ given by (\ref{defirrA}).
 \qed

 \bibliographystyle{plain}

\begin{thebibliography}{MMM}
 \bibitem{BoH91} M. Boyle and D. Handelman, The spectra of nonnegative matrices via symbolic dynamics,
 \emph{Annals of Math.} 133 (1991), 249--316.
 \bibitem{Fri78} S. Friedland,  On inverse problem for nonnegative and eventually
 nonnegative matrices, \emph{Israel J. Math.} 29 (1978), 43-60.
 \bibitem{Gan59} F.R. Gantmacher, {\it The Theory of Matrices}, Vol. I and II, Chelsea Publ. Co.,
 New York 1959.
 \bibitem{HoJ88} R.A. Horn and C.R. Johnson, {\it Matrix Analysis},
 Cambridge Univ. Press, New York 1988.
 \bibitem{KOR00} K.H. Kim, N.S. Ormes and F.W. Roush, The spectra of nonnegative integer matrices via formal power series,
 \emph{J. Amer. Math. Soc.}  13  (2000), 773--806.
 \bibitem{Laf10} T.J. Laffey, A constructive version of the Boyle-Handelman theorem on the spectra of nonnegative matrices,
 arXiv:1005.0929.  
 \bibitem{LoL78} R. Loewy and D. London, A note on an inverse problem for nonnegative matrices,
 \emph{Linear Multilin. Algebra} 6 (1978) 83--90.
 \bibitem{LaM98} T. Laffey and E. Meehan, A refinement of an inequality of Johnson, Loewy and London on nonnegative matrices and some applications,
 \emph{Electron. J. Linear Algebra}  3  (1998), 119--128
 \bibitem{TM07} J. Torre-Mayo, J., M.R. Abril-Raymundo, E. Alarcia-Estévez, C. Marijuán and M. Pisonero,
 The nonnegative inverse eigenvalue problem from the coefficients of the characteristic polynomial. EBL digraphs,
 \emph{Linear Algebra Appl.} 426 (2007), no. 2-3, 729--773.












 \end{thebibliography}

\end{document}